\documentclass[12pt]{amsart}
\usepackage{amssymb,mathrsfs}
\usepackage{amsthm,amsmath,amssymb}
\usepackage[numbers,sort&compress]{natbib}

\makeatletter

\newcommand{\Rmnum}[1]{\expandafter\@slowromancap\romannumeral #1@}
\makeatother
\usepackage{mathtools}

\theoremstyle{plain}
\newtheorem{theorem}{Theorem}[section]
\newtheorem{lemma}[theorem]{Lemma}
\newtheorem{proposition}[theorem]{Proposition}
\newtheorem{corollary}[theorem]{Corollary}

\theoremstyle{definition}

\newcommand{\abs}[1]{\lvert#1\rvert}

\usepackage[dvips]{graphics}
\usepackage{mathrsfs}
\usepackage{amssymb}
\usepackage{stmaryrd}
\usepackage{pifont}
\usepackage{amsmath}
\usepackage{hyperref}
\usepackage{amsthm}
\usepackage{cases}
\usepackage{mathtools}
\allowdisplaybreaks 
\title[unbounded continuous operators]{Continuous operators for unbounded convergence in Banach lattices}

\date{\today}


\author[Z. Wang]{Zhangjun Wang$^{1}$}
\address{$^1$ The first author:School of Mathematics, Southwest Jiaotong University,
	Chengdu, Sichuan,
	China, 611756.}
\email{zhangjunwang@my.swjtu.edu.cn}

\author[Z. Chen]{Zili Chen$^{2}$}
\address{$^2$ The second author:School of Mathematics, Southwest Jiaotong University, Chengdu, Sichuan,
	China, 611756.}
\email{zlchen@home.swjtu.edu.cn}

\author[J. Chen]{Jinxi Chen$^{3}$}
\address{$^3$ The third author: School of Mathematics, Southwest Jiaotong University, Chengdu, Sichuan,
	China, 611756.}
\email{jinxichen@home.swjtu.edu.cn}

\begin{document}

\begin{abstract}
Recently, the functionals different types of unbounded convergences $(uo, un, u\\aw, uaw^*)$ in Banach lattices were studied. In this paper, we study the continuous operators with respect to unbounded convergences. We first investigate the approximation property of continuous operators for unbounded convergence. Then we show some characterizations of the continuity of the continuous operators for $uo$, $un$, $uaw$ and $uaw^*$-convergence. Based on these results, we discuss the order-weakly compact operators on Banach lattices. Some related results are obtained as well.
\end{abstract}
	
\maketitle

\section{Introduction}
A net $(x_\alpha)$ in a Banach lattice $E$ is \emph{unbounded order (resp. norm, absolute weak)} convergent to some $x$, denoted by $x_\alpha\xrightarrow{uo}x$ (resp. $x_\alpha\xrightarrow{un}x$, $x_\alpha\xrightarrow{uaw}x$), if the net $(|x_\alpha-x|\wedge u)$ converges to zero in order (resp. norm, weak) for all $u\in E_+$. A net $(x_\alpha^\prime)$ in a dual Banach lattice $E^\prime$ is unbounded absolute weak* convergent to some $x^\prime$, denoted by $x_\alpha^\prime\xrightarrow{uaw^*}x^\prime$, if $|x^\prime_\alpha-x^\prime|\wedge u^\prime\xrightarrow{w^*}0$ for all $u^\prime\in E_+^\prime$. 
For the basic theory of $uo$, $un$, $uaw$ and $uaw^*$-convergence, we refer to \cite{GTX:16,KMT:16,Z:16,T:18}.

It can be easily verified that, in $l_p (1\leq p<\infty)$, $uo$, $un$ and $uaw$ and $uaw^*$-convergence of nets are the same as coordinate-wise convergence.  In $L_p(\mu) (1\leq p<\infty)$ for finite measure $\mu$, $uo$-convergence for sequences is the same as almost everywhere convergence, $un$ and $uaw$-convergence for sequences are the same as convergence in measure. In $L_p(\mu) (1< p<\infty)$ for finite measure $\mu$, $uaw^*$-convergence for sequences is also the same as convergence in measure.

In \cite{W:21}, we studied the continuity of the linear functionals for different types of unbounded convergences $(uo; un; uaw; uaw^*)$ in Banach lattices. The aims of the present paper are the continuity of linear operators for different types of unbounded convergences $(uo; un; uaw; uaw^*)$ in Banach lattices. In \cite[Theorem~2.7]{W:19}, we dicussed the "unbounded-norm" continuous operators and had the following. 
\begin{theorem}\label{M-weakly compact operators}
(\cite[Theorem~2.7]{W:19}) Let $E$ and $F$ be Banach lattices, for a continuous operator $T:E\rightarrow F$, the following statements hold.
	\begin{enumerate}
		\item The following conditions are equivalent.
		\begin{enumerate}
			\item $T$ is a M-weakly compact.
			\item $Tx_n\rightarrow0$ for every $uaw$-null sequence $(x_n)$ in $B_E$.
			\item $Tx_n\rightarrow0$ for every $uo$-null sequence $(x_n)$ in $B_E$.
		\end{enumerate}  
		\item For its adjoint operator $T^\prime:F^\prime\rightarrow E^\prime$, the following conditions are equivalent.
		\begin{enumerate}
			\item $T^\prime:F^\prime\rightarrow E^\prime$ is a M-weakly compact operator.
			\item $T^\prime y^\prime_n\rightarrow0$ for every $uaw^*$-null sequence $(y_n^\prime)$ in $B_{F^\prime}$.		
			\item $T^\prime y^\prime_n\rightarrow0$ for every $uaw$-null sequence $(y_n^\prime)$ in $B_{F^\prime}$.
			\item $T^\prime y^\prime_n\rightarrow0$ for every $uo$-null sequence $(y_n^\prime)$ in $B_{F^\prime}$.	
		\end{enumerate}
	\end{enumerate}
	
\end{theorem}
Now we study the "unbounded-unbounded" continuous operators. In the first part of the paper, we investigate the approximation property of continuous operators for unbounded convergence. Then we show some characterizations of the continuity of the continuous operators for $uo$, $un$, $uaw$ and $uaw^*$-convergence. As an application of these results, we conclude the paper with characterizations of the order-weakly compact operators on Banach lattices. Some related results are obtained as well.

Recall that a Riesz space $E$ is an ordered vector space in which $x\vee y=\sup\{x,y\}$ and $x\wedge y=\inf\{x,y\}$ exists for every $x,y\in E$. The positive cone of $E$ is denoted by $E_+$, $i.e.,E_+=\{x\in E:x\geq0\}$. For any vector $x$ in $E$ define $x^+:=x\vee0,x^-:=(-x)\vee0, |x|:=x\vee(-x)$. An operator $T:E\rightarrow F$ between two Riesz spaces is said to be \emph{positive} if $Tx\geq0$ for all $x\geq0$. A net $(x_\alpha)$ in a Riesz space is called disjoint whenever $\alpha\ne \beta$ implies $|x_\alpha|\wedge|x_\beta|=0$ (denoted by $x_\alpha\perp x_\beta$). A set $A$ in $E$ is said to be \emph{order bounded} if there exists some $u\in E_+$ such that $|x|\leq u$ for all $x\in A$. The solid hull $Sol(A)$ of $A$ is the smallest solid set including $A$ and it equals the set $Sol(A):=\{x\in E:\exists y\in A,|x|\leq|y|\}.$ An operator $T:E\rightarrow F$ is called order bounded if it maps order bounded subsets of $E$ to order bounded subsets of $F$. A Banach lattice $E$ is a Banach space $(E, \Vert\cdot\Vert)$ such that $E$ is a Riesz space and its norm satisfies the following property: for each $x, y\in E$ with $|x|\leq|y|$, we have $\Vert x\Vert\leq\Vert y\Vert$.

For undefined terminology, notation and basic theory of Riesz space, Banach lattice and linear operator, we refer to \cite{AB:06,MN:91}.

\section{Results}\label{}
A \emph{Riesz pseudonorm} is a real-valued function defined on a Riesz space $E$ ($\rho:E\rightarrow \mathbb{R}$) satisfying the following properties.
\begin{enumerate}
\item $\rho(x)\geq0$ for all $x\in E$.
\item $\rho(x+y)\leq\rho(x)+\rho(y)$ for all $x,y\in E$.
\item $\rho(\lambda_n x)\rightarrow0$ as $\lambda_n\rightarrow0$ for all $x\in E$.
\item $\rho(x)\leq\rho(y)$ whenever $\abs{x}\leq\abs{y}$ holds in $E$.
\end{enumerate}

A linear topology on a Riesz space is locally solid iff it is generated by a family of Riesz pseudonorms. It is natural to consider the "unbounded" topology generated by a family of Riesz pseudonorms.

\begin{proposition}\label{un,uaw,uaw^*-topology}
Let $(E,\Vert\cdot\Vert)$ be a Banach lattice, the following statements hold.
\begin{enumerate}
\item The map $\rho_u:E\rightarrow \mathbb{R_+}$ defined by $\rho_u(x)=\big\Vert\abs{x}\wedge u\big\Vert$ is a Riesz pseudonorm for each $u\in E_+$. Moreover, the $un$-topology and the topology generated by $(\rho_u)_{u\in E_+}$ coincide.
\item The map $\rho_{u,x^\prime}:E\rightarrow R_+$ as $\rho_{u,x^\prime}(x)=x^\prime(\abs{x}\wedge u)$ is a Riesz pseudonorm for each $u\in E_+,x^\prime\in E^\prime_+$. Moreover, the $uaw$-topology and the topology generated by the family $(\rho_{u,f})_{u\in E_+,x^\prime\in E^\prime_+}$ coincide.
\item The map $\rho_{u^\prime,x}:E^\prime\rightarrow R_+$ as $\rho_{u^\prime,x}(x^\prime)=(\abs{x^\prime}\wedge u^\prime)(x)$ is a Riesz pseudonorm for each $u^\prime\in E^\prime_+,x\in E_+$. Moreover, the $uaw^*$-topology and the topology generated by the family $(\rho_{u^\prime,x})_{u^\prime\in E^\prime_+,x\in E_+}$ coincide.
\end{enumerate} 
\end{proposition}	
\begin{proof}
We just prove $(3)$, the rest of proof is similar.

Clearly, $\rho_{u^\prime,x}(x^\prime)\geq0$ for all $x^\prime\in E^\prime$ and $\rho_{u^\prime,x}(x^\prime)\leq\rho_{u^\prime,x}(y^\prime)$ whenever $\abs{x^\prime}\leq\abs{y^\prime}$ holds in $E^\prime$. Since $\abs{x^\prime+y^\prime}\wedge u^\prime\leq (\abs{x^\prime}+\abs{y^\prime})\wedge u^\prime\leq\abs{x^\prime}\wedge u^\prime+\abs{y^\prime}\wedge u^\prime$, we have $\rho_{u^\prime,x}(x^\prime+y^\prime)\leq\rho_{u^\prime,x}(x^\prime)+\rho_{u^\prime,x}(y^\prime)$. According to the inequality $\rho_{u^\prime,x}(\lambda_n x^\prime)=(\abs{\lambda_n x^\prime}\wedge u^\prime)(x)\leq(\abs{\lambda_n x^\prime})(x)=\abs{\lambda_n}|x^\prime|(x)$, we have $\rho_{u^\prime,x}(\lambda_n x^\prime)\rightarrow0$ whenever $\lambda_n\rightarrow0$. Therefore, $\rho_{u,x^\prime}$ is a Riesz pseudonorm for all $u\in E_+,x^\prime\in E^\prime_+$. It is easy to see that the $uaw^*$-topology and the topology generated by the family $(\rho_{u^\prime,x})_{u^\prime\in E^\prime_+,x\in E_+}$ coincide.
\end{proof}

The following theorems describe an approximation property of continuous operators for unbounded convergence.
\begin{theorem}\label{disjoint-u}
Let $T:E\rightarrow F$ be a continuous operator from a Banach lattice $E$ to a Banach lattice $F$, let $A$ be a bounded solid subset of $E$, and let $(\rho_\gamma)_{\gamma\in \Gamma}$ be a family of norm continuous Riesz pseudonorms on $F$. If $\rho_\gamma(Tx_n)\rightarrow0$ holds for every disjoint sequence $(x_n)\subset A$ and all $\gamma\in \Gamma$, then for each $\epsilon>0,\gamma\in \Gamma$ there exists some $u\in E_+$ lying the ideal generated by $A$ such that $\rho_\gamma\big(T(\abs{x}-u)^+\big)<\epsilon$ holds for all $x\in A$.
\end{theorem}	
\begin{proof}
If the claim is not true, then there exists some $\epsilon>0,\gamma\in \Gamma$ such that, for any $u\geq0$ in the ideal generated by $A$, we have $\rho_\gamma\big(T(\abs{x}-u)^+\big)\geq\epsilon$ for at least one $x\in A$. In particular, there exists a sequence $(x_n)\subset A$ such that, for all $n\in\mathbb{N}$, we have $\rho_\gamma\big(T(\abs{x_{n+1}}-4^n\sum_{i=1}^{n}\abs{x_i})^+\big)\geq\epsilon$.
	
Put $y=\sum_{n=1}^{\infty}2^{-n}\abs{x_n}$, $w_n=(\abs{x_{n+1}}-4^n\sum_{i=1}^{n}\abs{x_i})^+$ and $v_n=(\abs{x_{n+1}}-4^n\sum_{i=1}^{n}\abs{x_i}-2^{-n}y)^+$. According to \cite[Lemma~4.35]{AB:06}, $(v_n)$ is disjoint. $(v_n)\subset A$ since $A$ is solid and $0\leq v_n\leq \abs{x_{n+1}}$. By our hypothesis, we have $\rho_\gamma(Tv_n)\rightarrow0$. It follows form $0\leq w_n-v_n\leq2^{-n}y$ that $\Vert w_n-v_n\Vert\leq2^{-n}\Vert y\Vert$. Since $\rho_\gamma$ is a norm continuous Riesz pseudonorm, hence $\rho_\gamma(T(w_n-v_n))\rightarrow0$. From $\rho_\gamma(Tw_n)\leq\rho_\gamma(T(w_n-v_n))+\rho_\gamma(Tv_n)$, we see that $\rho_\gamma(Tw_n)=\rho_\gamma\big(T(\abs{x_{n+1}}-4^n\sum_{i=1}^{n}\abs{x_i})^+\big)\rightarrow0$. This leads to a contradiction, and the proof is completed. 
\end{proof}

\begin{corollary}\label{disjoint-un,uaw,uaw^*}
Let $T:E\rightarrow F$ be a continuous operator from a Banach lattice $E$ to a Banach lattice $F$ and $A$ a bounded solid subset of $E$, the following statements hold.
\begin{enumerate}
\item\label{disjoint-un} If $Tx_n\xrightarrow{un}0$ (i.e. $\rho_v(Tx_n)=\Vert\abs{Tx_n}\wedge v\Vert\rightarrow0$ for all $v\in F_+$) holds for every disjoint sequence $(x_n)\subset A$, then for each $\epsilon>0,v\in F_+$ there exists some $u\in E_+$ lying the ideal generated by $A$ such that $\rho_v\big(T(\abs{x}-u)^+\big)<\epsilon$ holds for all $x\in A$.
\item\label{disjoint-uaw} If $Tx_n\xrightarrow{uaw}0$ (i.e. $\rho_{v,y^\prime}(Tx_n)=y^\prime(\abs{Tx_n}\wedge v)\rightarrow0$ for all $v\in F_+,y^\prime\in F^\prime_+$) holds for every disjoint sequence $(x_n)\subset A$, then for each $\epsilon>0,v\in F_+,y^\prime\in F^\prime_+$ there exists some $u\in E_+$ lying the ideal generated by $A$ such that $\rho_{v,y^\prime}\big(T(\abs{x}-u)^+\big)<\epsilon$ holds for all $x\in A$.
\item\label{disjoint-uaw^*} If $F$ is a dual Banach lattice (denoted by $F^\prime$ in here) and $Tx_n\xrightarrow{uaw^*}0$ (i.e. $\rho_{v^\prime,y}(Tx_n)=(\abs{Tx_n}\wedge v^\prime)(y)\rightarrow0$ for all $v^\prime\in F^\prime_+,y\in F_+$) holds for every disjoint sequence $(x_n)\subset A$, then for each $\epsilon>0,v^\prime\in F^\prime_+,y\in F_+$ there exists some $u\in E_+$ lying the ideal generated by $A$ such that $\rho_{v^\prime,y}\big(T(\abs{x}-u)^+\big)<\epsilon$ holds for all $x\in A$.
\end{enumerate}
\end{corollary}		
The following results show some characterizations of the continuity of the continuous operators for $uo$, $un$, $uaw$ and $uaw^*$-convergence.

Recall that an operator $T : E\rightarrow F$ from a Riesz space $E$ to a Banach space $F$ is called \emph{order-weakly compact} if $T [-x , x]$ is relatively weakly compact for all $x\in E_+$ . A continuous operator $T : E\rightarrow F$ between two Banach lattices is said to be \emph{preserve a sublattice isomorphic to $L$} if there exists closed sublattices $U\subset E , V\subset F$ isomorphic to $L$ such that the restriction of $T$ to $U$ acts as an isomorphism onto $V$.
\begin{lemma}\label{d-un is ow}
For a continuous operator $T:E\rightarrow F$ from a Dedekind $\sigma$-complete Banach lattice $E$ to a Banach lattice $F$, if $Tx_n\xrightarrow{un}0$ for every disjoint sequence $(x_n)\subset B_E$, then $T$ is order-weakly compact.
\end{lemma}	
\begin{proof}
Assume that $T:E\rightarrow F$ is not order-weakly compact, according to \cite[Corollary~3.4.5]{MN:91}, $T$ preserve a sublattice isomorphic to $\ell_\infty$. Let $(e_n)$ be the unit vectors of $\ell_\infty$, clearly, $(Te_n)$ is not $un$-null in $\ell_\infty$. Indeed, if $(Te_n)$ is $un$-null in $\ell_\infty$, then $Te_n\rightarrow0$ in $\ell_\infty$ by \cite[Theorem~2.3]{KMT:16}, moreover $e_n\rightarrow0$ since the restriction of $T$ to $\ell_\infty$ is isomorphic. Hence, $(Te_n)$ is not $un$-null in $F$. This leads to a contradiction. Therefore, $T$ is order weakly compact.
\end{proof}
\begin{theorem}\label{u-u=d-u}
Let $E$ and $F$ be Banach lattices, for a continuous operator $T:E\rightarrow F$, the following statements hold.
\begin{enumerate}
\item If $E$ is Dedekind $\sigma$-complete, then the following conditions are equivalent.
\begin{enumerate}
\item $Tx_n\xrightarrow{un}0$ for every $uo$-null sequence $(x_n)\subset B_E$.
\item $Tx_n\xrightarrow{un}0$ for every $uaw$-null sequence $(x_n)\subset B_E$.
\item $Tx_n\xrightarrow{un}0$ for every disjoint sequence $(x_n)\subset B_E$.
\end{enumerate}
\item If $T$ is positive, then the following conditions are equivalent.
\begin{enumerate}
\item $Tx_n\xrightarrow{uaw}0$ for every $uaw$-null sequence $(x_n)\subset B_E$.
\item $Tx_n\xrightarrow{uaw}0$ for every disjoint sequence $(x_n)\subset B_E$.
\end{enumerate}
\item For the adjoint $T^\prime:F^\prime\rightarrow E^\prime$ of positive operator $T$, the following conditions are equivalent.
\begin{enumerate}
\item $T^\prime y_n^\prime\xrightarrow{uaw^*}0$ for every $uaw^*$-null sequence $(y_n^\prime)\subset B_{F^\prime}$.
\item $T^\prime y_n^\prime\xrightarrow{uaw^*}0$ for every $uaw$-null sequence $(y_n^\prime)\subset B_{F^\prime}$.
\item $T^\prime y_n^\prime\xrightarrow{uaw*}0$ for every $uo$-null sequence $(y_n^\prime)\subset B_{F^\prime}$.
\item $T^\prime y_n^\prime\xrightarrow{uaw^*}0$ for every disjoint sequence $(y_n^\prime)\subset B_{F^\prime}$.
\end{enumerate}
\end{enumerate}
\end{theorem}
\begin{proof}
$(1)(a)\Rightarrow(1)(c)$ and $(1)(b)\Rightarrow(1)(c)$. According to \cite[Corollary~3.6]{GTX:16} and \cite[Lemma~2]{Z:16}, every disjoint sequence is $uo$-null and $uaw$-null.
	
$(1)(c)\Rightarrow(1)(a)$ and $(1)(c)\Rightarrow(1)(b)$. Assume that $(3)$ holds, $T$ is order-weakly compact by Lemma \ref{d-un is ow}. It follows from Corollary \ref{disjoint-un,uaw,uaw^*} that, for each $\epsilon>0,v\in F_+$, there exists some $u\in E_+$ such that $\rho_v\big(T(\abs{x}-u)^+\big)=\rho_v\big(T(\abs{x})-T(\abs{x}\wedge u)\big)<\epsilon$ holds for all $x\in B_E$.	

For a $uo$-null (resp. $uaw$-null) sequence $(x_n)\subset B_E$, clearly, $\abs{x_n}\wedge u\xrightarrow{o (resp. w)}0$. Since $T$ is order-weakly compact, hence                                                                                                                                                                                                                                                                                                                                                                                                                                                                                                                                                                                                    $T(\abs{x_n}\wedge u)\rightarrow0$ by \cite[Theorem~3.4.4 and Corollary~3.4.9]{MN:91}, moreover $\rho_v\big(T(\abs{x_n}\wedge u)\big)\rightarrow0$. Thus, $T\abs{x_n}\xrightarrow{un}0$.	Since $x_n\xrightarrow{uo (resp. uaw)}$ iff $x_n^\pm\xrightarrow{uo (resp. uaw)}0$ iff $\abs{x_n}\xrightarrow{uo (resp. uaw)}0$, therefore $Tx_n\xrightarrow{un}0$.

$(2)(a)\Rightarrow(2)(b)$. Obvious.

$(2)(b)\Rightarrow(2)(a)$. It is similar to the proof os $(1)(c)\Rightarrow(1)(b)$ that, for each $\epsilon>0,v\in F_+,y^\prime\in F^\prime_+$ there exists some $u\in E_+$ such that $\rho_{v,y^\prime}\big(T(\abs{x}-u)^+\big)<\epsilon$ holds for all $x\in B_E$. Using $T$ is weak-weak continuous and positive, the rest of proof is similar.

$(3)(a)\Rightarrow(3)(b)$ and $(3)(a)\Rightarrow(3)(c)$. Every $uo$ (resp. $uaw$-null) sequence is $uaw^*$-null.

$(3)(b)\Rightarrow(3)(d)$ and $(3)(c)\Rightarrow(3)(d)$. Clearly.

$(3)(d)\Rightarrow(3)(a)$. According $T^\prime$ is positive and weak*-weak* continuous, the proof is similar.
\end{proof}

Recall that a closed sublattice $L$ in Banach lattice $E$ which is generated by disjoint sequence $(x_n)\subset E_+$ is called \emph{sup-stable} if $x=\sum_{n=1}^{\infty}a_n x_n$ for every $x\in E$ with $\lim\sup\abs{a_n}<\infty$. For the unit vectors $(e_n)$ of $\ell_\infty$, clearly, $c_0$ is a closed sublattice in $\ell_\infty$ generated by $\big(e_n\wedge(1,\frac{1}{2},...,\frac{1}{n},...)\big)_{n=1}^\infty$ and $\ell_\infty$ is a closed sup-stable sublattice generated by $\big(e_n\wedge(1,1,1,...)\big)_{n=1}^\infty$. It is natural to ask that whether the result of the example holds in more general situations. The following results confirm the hypothesis.
\begin{proposition}\label{iso-cl}
Suppose that $E$ is a Banach lattice and $(x_n)\subset B_E$ is a disjoint sequence but not $un$-null, the following statements hold.
\begin{enumerate}
\item The closed sublattice $L$ generated by $\{\abs{x_n}\wedge u:n\in \mathbb{N}\}$ is isomorphic to $c_0$ for some $u\in E_+$.
\item If $E$ is Dedeking $\sigma$-complete, then the closed sup-stable sublattice $L$ generated by $\{\abs{x_n}\wedge u:n\in \mathbb{N}\}$ is isomorphic to $\ell_\infty$ for some $u\in E_+$.
\end{enumerate}
\end{proposition}
\begin{proof}
Since $(x_n)$ is not $un$-null, hence there exists some $\epsilon>0$ and a Riesz pseudonorm $\rho_u$ such that $\rho_u(x_n)=\big\Vert\abs{x_n}\wedge u\big\Vert>\epsilon$ for all $n$. Without loss of generality we may assume that $\rho_u(x_n)=\big\Vert\abs{x_n}\wedge u\big\Vert=1$.
	
$(1)$ For $n\in\mathbb{N}$ and $z=(a_i)_{i=1}^n\in\mathbb{R}^n$. Since $\Vert z\Vert_\infty\leq\Vert\sum_{i=1}^{n}a_i\cdot (\abs{x_n}\wedge u)\Vert\leq\Vert u\Vert\cdot\Vert z\Vert_\infty$, therefore the closed sublattice $L$ generated by $\{\abs{x_n}\wedge u\}$ is isomorphic to $c_0$.
	
$(2)$ For any $0\leq z=(a_i)_{i=1}^\infty\in \ell_\infty$, let $j(z)=\sup\{a_n\cdot(\abs{x_n}\wedge u)\}\in L$, the supremum exsits since $E$ is Dedekind $\sigma$-complete. Clearly, $j:\ell^\infty_+\rightarrow L_+$ is additive and positively homogeneous. Thus, $j$ extends to all of $\ell_\infty$ as a lattice and norm isomorphism.
\end{proof}
An operator $T:E\rightarrow F$ is called \emph{unbounded norm continuous} ($un$-continuous, for short) if $Tx_n\xrightarrow{un}0$ for every $un$-null sequence $(x_n)\subset B_E$. According to $\abs{Tx_n}\wedge v=T\abs{x_n}\wedge Tu=T(\abs{x_n}\wedge u)$ for all $u\in E_+, v\in F_+$, $x_n\xrightarrow{un}0$ implies $Tx_n\xrightarrow{un}0$, therefore every onto lattice homomorphism is $un$-continuous. For every measurable space $(\Omega,\varSigma)$, let $St(\varSigma)$ denote the collection of all measurable stepfuntions: $f=\sum_{i=1}^{r}a_i\chi_{A_i}:\Omega\rightarrow \mathbb{R}$ where $r\in\mathbb{N}$, $a_i\in\mathbb{R}$ and $A_i\in \varSigma$ for $i\in\{1,2,...r\}$.
\begin{theorem}\label{T-iso-cl}
Let $E$ and $F$ be Banach lattices, for a onto lattice homomorphism $T:E\rightarrow F$. Assume that $(Tx_n)$ is not $un$-null in $F$ for a disjoint sequence $(x_n)\subset B_E$, then there exists a subsequence $(k(n))_1^\infty$ such that the following statements hold.
\begin{enumerate}
\item The closed sublattice $L$ generated by $\{\abs{x_{k(n)}}\wedge u:n\in \mathbb{N}\}$ is isomorphic to $c_0$ and $T$ acts on $L$ as a isomorphism for some $u\in E_+$.
\item If $E$ is Dedekind $\sigma$-complete, then the closed sup-stable sublattice $L$ generated by $\{\abs{x_{k(n)}}\wedge u:n\in \mathbb{N}\}$ is isomorphic to $l_\infty$ and $T$ acts on $L$ as a isomorphism for some $u\in E_+$.
\end{enumerate}
\end{theorem}
\begin{proof}
Since $E^{\prime\prime}$ is Dedekind complete and $E$ is a sublattice of $E^{\prime\prime}$, hence $(2)$ implies $(1)$. Therefore, we just prove $(2)$.
	
Since $(Tx_n)$ is not $un$-null, so there exists a Riesz pseudonorm $\rho_v$ such that $\rho_v(Tx_n)=\big\Vert\abs{Tx_n}\wedge v\big\Vert>1$ for all $n$ and some $v\in F_+$. Since $T$ is onto and lattice homomorphism, therefore $(x_n)$ is not $un$-null, that is, $\abs{x_n}\wedge u\nrightarrow0$ for some $u\in E_+$ satisfying $Tu=v$. For every $n\in N$, let $y_n^\prime\in B_{F^\prime}\cap F^\prime_+$ such that $y_n^\prime(\abs{Tx_n}\wedge v)>1$. We set $x^\prime_n=T^\prime y_n^\prime$, since $x_n^\prime(\abs{x_n}\wedge u)=y_n^\prime(\abs{Tx_n}\wedge v)$, so $x_n^\prime(\abs{x_n}\wedge u)>1$.
	
Let $0<\epsilon<\frac{1}{4}$. Cleary, $(\abs{x_n}\wedge u)_1^\infty\subset[0,u]\subset E_+$, $\sup_n\big|\abs{x_n^\prime}(u)\big|<\infty$ and $\sum_{1}^{\infty}\abs{x_n}\wedge u\leq u$. It follows from \cite[Theorem~2.3.7]{MN:91} that there exists a subsequence $(k(n))_1^\infty$ satisfying $\big|\abs{x^\prime_{k(n)}}(\sup\{\abs{x_{k(j)}}\wedge u:n\ne j\in \mathbb{N}\})\big|<\epsilon$. We may assume that $k(n)=n$ for all $n$. For every finitely valued $z\in l_\infty$, there exists disjoint subsets $A_1,A_2,...,A_r\subset N$ and $a_1,a_2,...,a_r\in R$ such that $z=\sum_{i=1}^{r}a_i\cdot \chi_{A_i}\in St(\mathbb{N}                                                                                                                                         )$. For all $i=1,2,...,r$, we define $v_i=\sup\{\abs{x_n}\wedge u:n\in A_i\}$ and $T_1z=\sum_{i=1}^{r}a_i Tv_i$. 

Clearly, $T_1$ is well defined and linear on $St(\mathbb{N})$ such that $\abs{y^\prime(T_1z)}\leq\sum_{i=1}^{r}\abs{a_i}\cdot\abs{y^\prime(Tv_i)}\leq\Vert z\Vert_\infty\Vert T\Vert\Vert u\Vert$ for all $y^\prime\in B_{F^\prime}$. Therefore, $\Vert T_1\Vert\leq \Vert T\Vert\Vert u\Vert$.
	
We may assume that $\Vert z\Vert_\infty=\abs{a_1}$. For every $n\in A_1$,
\begin{align*}
\Vert T_1z\Vert&\geq\abs{y_n^\prime(T_1z)}=\abs{x^\prime_n(\sum_{i=1}^{r}a_i\cdot v_i)}\\
&\geq\Vert z\Vert_\infty\cdot\abs{x^\prime_n(\abs{x_n}\wedge u)}-\Vert z\Vert_\infty\cdot\big|\abs{x^\prime_n}(\sup\{\abs{x_j}\wedge u:j\ne n\})\big|\\
&\geq\Vert z\Vert_\infty(1-\epsilon)\geq\dfrac{\Vert z\Vert_\infty}{2}.
\end{align*}
	
Since $St(\mathbb{N})$ is dense in $\ell_\infty$, hence $T_1$ extends as an lattice and norm isomorphism to all of $l_\infty$. According to Theorem \ref{iso-cl}, the sup-stable sublattice generated by $(\abs{x_n}\wedge u)$ is isomorphic to $l_\infty$. The proof is completed.
\end{proof}
\begin{theorem}\label{u-un=d-un=o-w}
Let $E$ be a Dedekind $\sigma$-complete Banach lattice and $F$ a Banach lattice, for a onto lattice homomorphism $T:E\rightarrow F$, the following conditions are equivalent.
\begin{enumerate}
\item $Tx_n\xrightarrow{un}0$ for every $uo$-null sequence $(x_n)\subset B_E$.
\item $Tx_n\xrightarrow{un}0$ for every $uaw$-null sequence $(x_n)\subset B_E$.
\item $Tx_n\xrightarrow{un}0$ for every disjoint sequence $(x_n)\subset B_E$.
\item $T$ is order-weakly compact.
\end{enumerate}
\end{theorem}
\begin{proof}
$(1)\Leftrightarrow(2)\Leftrightarrow(3)$ by Theorem \ref{u-u=d-u}.
	
$(3)\Rightarrow(4)$ by Lemma \ref{d-un is ow}.
	
$(4)\Rightarrow(3)$. Assume that $(3)$ does not hold, then there exists a bounded disjoint sequence $(x_n)$ such that $(Tx_n)$ is not $un$-null. According to Theorem \ref{T-iso-cl}, $T$ preserves a sublattice isomorphic to $l_\infty$. It follows from \cite[Corollary~3.4.5]{MN:91} that $T$ is not order-weakly compact. This leads to a contradiction.

$(4)\Rightarrow(1)$ and $(4)\Rightarrow(2)$ by $\abs{Tx_n}\wedge v=T\abs{x_n}\wedge Tu=T(\abs{x_n}\wedge u)$ for all $u\in E_+, v\in F_+$ and \cite[Theorem~3.4.4 and Corollary~3.4.9]{MN:91}. 
\end{proof}
\noindent \textbf{Acknowledgement.} The research is supported by National Natural Science Foundation of China(NSFC:51875483).

\end{document}